\documentclass{amsart}
\usepackage{amsmath,amssymb,graphicx,psfrag}

\newcommand{\wP}{\widetilde{P}}
\newcommand{\wS}{\widetilde{S}}

\newcommand\Z{\mathbb{Z}}
\newcommand\R{\mathbb{R}}

\newcommand\cA{{\mathcal A}}
\newcommand\cF{{\mathcal F}}
\newcommand\cS{{\mathcal S}}
\newcommand\cC{{\mathcal C}}
\newcommand\fF{\mathfrak F}

\newcommand{\A}[1]{a^{(#1)}}
\newcommand{\B}[1]{b^{(#1)}}

\newcommand{\tropical}[1]{\left[\,#1\,\right]}

\newcommand{\I}{^{-1}}

 \newcommand{\pichere}[2]
 {\begin{center}\includegraphics[width=#1\textwidth]{#2}\end{center}}

 \newcommand{\lab}[3]{\psfrag{#1}[#3]{$\scriptstyle{#2}$}}

\newtheorem{thm}{Theorem}
\newtheorem{lem}[thm]{Lemma}

\theoremstyle{definition}

\theoremstyle{remark}

\newtheorem{remark}[thm]{Remark}

\allowdisplaybreaks[1]

\graphicspath{{./}{figures/}}

\begin{document}

\title{On the topological entropy of families of braids}
\author{Toby Hall} 
\email{tobyhall@liv.ac.uk}
\address{Department of Mathematical Sciences\\ University of Liverpool\\
  Liverpool L69 7ZL\\ UK}

\author{S. \"Oyk\"u Yurtta\c s}
\email{oyku1981@liv.ac.uk}
\address{Department of Mathematical Sciences\\ University of Liverpool\\
  Liverpool L69 7ZL\\ UK}
\thanks{S.O.Y. is grateful to Dicle University for sponsorship}

\begin{abstract}
A method for computing the topological entropy of each braid in an
infinite family, making use of Dynnikov's coordinates on the boundary
of Teichm\"uller space, is described. The method is illustrated on two
two-parameter families of braids.
\end{abstract}
\keywords{Topological entropy, Nielsen-Thurston classification,
  Dynnikov coordinates}
\subjclass[2000]{37E30, 37B40}

\maketitle
\section{Introduction}
In the dynamical study of iterated surface homeomorphisms, it is
common to seek to compute the topological entropy of each member of an
infinite family of isotopy classes, perhaps on varying surfaces ---
the topological entropy of an isotopy class being the minimum
topological entropy of a homeomorphism in the class, which is realised
by a Nielsen-Thurston canonical
representative~\cite{Thurston,FLP,IsotStab}. The normal approach to
such a problem is to use train-track methods~\cite{BH,FM,Los}, which
not only make it possible to compute topological entropy, but also, in
the pseudo-Anosov case, provide a Markov partition for the
pseudo-Anosov homeomorphism in the isotopy class, and hence
information about the structure of its invariant singular measured
foliations.

One drawback of this approach is that even single train tracks are
fairly unwieldy objects. It is usually far from straightforward to
describe an infinite family of train tracks, to verify that they
are indeed invariant under the relevant isotopy classes, and to
compute the transition matrices and hence the topological entropy of
the induced train track maps: very often, the best that one can
reasonably do is to draw pictures of typical train tracks in the
family and rely on the reader's ability to observe that they are
invariant. 

In this paper an alternative approach to the problem is described in
the case of families of isotopy classes of orientation-preserving
homeomorphisms of punctured disks --- such isotopy classes can be
described by elements of Artin's braid groups. The method is
illustrated by applying it to two families of braids considered by
\mbox{Hironaka} and Kin~\cite{kin}, which are of interest in the study of
braids with low topological entropy. The results presented here about
these families are not new, therefore: the emphasis is on the method
used to obtain them, which can be contrasted with the train track
methods of Hironaka and Kin.

The methods developed are a relatively straightforward application of
Dynnikov's coordinate system~\cite{dynnikov} on the boundary of the
Teichm\"uller space of the punctured disk, together with the {\em
  update rules} which describe the action of the Artin braid
generators on the boundary of Teichm\"uller space in terms of Dynnikov
coordinates. This background material is described in
Section~\ref{sec:dynnikov}. The practical application of this theory
is very much eased by the results presented in
Section~\ref{sec:updatecontig}, which give update rules for braids
which can be written as ascending or descending sequences of
contiguous Artin generators (or their inverses), such as
$\sigma_3\sigma_4\sigma_5\sigma_6$. Examples of the application of the
method to two two-parameter families of braids are given in
Section~\ref{sec:main}: the examples include showing that an infinite
family of braids is of reducible type, as well as computing
topological entropies in the pseudo-Anosov case.

\section{Dynnikov coordinates of measured foliations}
\label{sec:dynnikov}
This section is essentially an expansion of parts of Dynnikov's very
terse paper~\cite{dynnikov}: see also~\cite{DDRW,efficient}, and~\cite{mous,FT}
for dynamical applications. One difference is that in the papers cited
above the action of the $n$-braid group $B_n$ on an $(n+2)$-punctured
disk is considered, whereas here, as is more appropriate in a
dynamical setting, $B_n$ acts on an $n$-punctured disk. This
modification requires separate consideration of the action of the
``end'' Artin generators $\sigma_1$ and $\sigma_{n-1}$. In addition,
the useful Lemma~\ref{lem:dynninvert} doesn't seem to have appeared
explicitly in the literature.

\subsection{The Dynnikov coordinates of a measured foliation}
Let~$D_n$ be a standard model of the $n$-punctured disk ($n\ge 3$).
Write~$\cF_n$ for the set of singular measured
foliations~$(\cF,\mu)$ on~$D_n$, and $\fF_n$ for~$\cF_n$ up to isotopy
and Whitehead equivalence (see for example~\cite{FLP}): the element of
$\fF_n$ containing $(\cF,\mu)\in \cF_n$ is
denoted~$[\cF,\mu]$. Dynnikov's coordinate system provides an explicit
bijection $\rho:\fF_n\to \R^{2n-4}\setminus\{0\}$.

Let~$(\cF,\mu)\in\cF_n$. Write $\cA_n$ for the set of arcs in~$D_n$
which have each endpoint either on the boundary or at a
puncture. Recall that if~$\alpha\in\cA_n$, then its measure
$\mu(\alpha)$ is defined to be
\[\mu(\alpha)=\sup \sum_{i=1}^k \mu(\alpha_i),\]
where the supremum is taken over all finite collections
$\alpha_1,\ldots, \alpha_k$ of mutually disjoint subarcs of~$\alpha$
which are transverse to~$\cF$. Denoting by $[\alpha]$ the isotopy
class of~$\alpha$ (under isotopies through~$\cA_n$), one can then define
\[\mu([\alpha])=\inf_{\beta\in[\alpha]}\mu(\beta),\]
which is well defined on~$\fF_n$.

\begin{figure}[htbp]
\lab{1}{1}{l}
\lab{2}{2}{l}
\lab{3}{i-1}{l}
\lab{4}{i}{l}
\lab{5}{i+1}{l}
\lab{6}{n-1}{l}
\lab{7}{n}{l}
\lab{a}{\alpha_1}{l}
\lab{b}{\alpha_2}{l}
\lab{c}{\alpha_{2i-5}}{r}
\lab{d}{\alpha_{2i-4}}{r}
\lab{e}{\alpha_{2i-3}}{l}
\lab{f}{\alpha_{2i-2}}{l}
\lab{g}{\alpha_{2i-1}}{l}
\lab{h}{\alpha_{2i}}{l}
\lab{i}{\alpha_{2n-5}}{r}
\lab{j}{\alpha_{2n-4}}{r}
\lab{k}{\beta_1}{tr}
\lab{l}{\beta_{n-1}}{tl}
\lab{m}{\beta_{i-1}}{br}
\lab{n}{\beta_i}{bl}
\pichere{0.95}{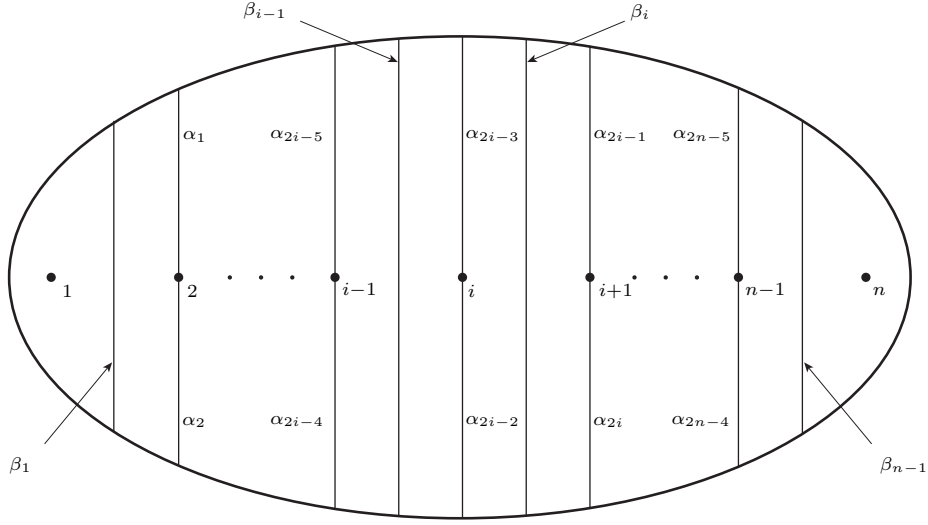}
\caption{The arcs $\alpha_i$ and $\beta_i$.}
\label{figure:munu}
\end{figure}

Consider the arcs $\alpha_i$ ($1\le i\le 2n-4$) and $\beta_i$ ($1\le
i\le n-1$) depicted in Figure~\ref{figure:munu}: the arcs
$\alpha_{2j-3}$ and $\alpha_{2j-2}$ (for $2\le j\le n-1$) join the
$j^{\text{th}}$ puncture to the boundary, while the arc $\beta_i$ has
both endpoints on the boundary and passes between the $i^{\text{th}}$
and $i+1^{\text{th}}$ punctures.

Let~$\tau:\fF_n\to\R_{\ge0}^{3n-5}$ be the {\em triangle coordinate function}
defined by 
\[\tau([\cF,\mu]) = \left(\mu([\alpha_1]), \ldots, \mu([\alpha_{2n-4}]),\,
\mu([\beta_1]), \ldots, \mu([\beta_{n-1}])\right).\] The
function~$\tau$ is injective: if $\tau([\cF,\mu])$ is given, then a
representative measured foliation in~$[\cF,\mu]$ can be constructed by
gluing together pieces of measured foliation in each of the strips of
Figure~\ref{figure:munu}. However, it is clearly not surjective:
$\tau([\cF,\mu])$ must satisfy the triangle inequality in each of the
strips of Figure~\ref{figure:munu}, as well as additional conditions
to ensure that~$(\cF,\mu)$ has no singularities which are centers.

Let~$\rho:\fF_n\to\R^{2n-4}\setminus\{0\}$ be the {\em Dynnikov
  coordinate function} defined by
\[\rho([\cF,\mu]) = (a,b)=(a_1,\ldots,a_{n-2},\,b_1,\ldots,b_{n-2}),\]
where for $1\le i\le n-2$
\[a_i=\frac{\mu([\alpha_{2i}])-\mu([\alpha_{2i-1}])}{2}\qquad
\text{and} \qquad b_i=\frac{\mu([\beta_i])-\mu([\beta_{i+1}])}{2}.\]
Let~$\cC_n=\R^{2n-4}\setminus\{0\}$ denote the space of Dynnikov
  coordinates. 

The Dynnikov coordinate function is a bijection (in fact it is a
homeomorphism when $\fF_n$ is endowed with its usual topology). To
describe its inverse, it is sufficient to describe a function
$\cC_n\to \R_{\ge0}^{3n-5}$ which sends each
$(a,b)\in\cC_n$ to the triangle coordinates of a
measured foliation~$[\cF,\mu]$ which has Dynnikov
coordinates~$(a,b)$. 

\begin{lem}[Inversion of Dynnikov coordinates]
\label{lem:dynninvert}
Let $(a,b)\in\cC_n$. Then $(a,b)$ is the Dynnikov
coordinate of exactly one element~$[\cF,\mu]$ of~$\fF_n$, which has
\begin{align*}
\mu([\beta_i]) &= 2\max_{1\le k\le n-2}\left(
|a_k|+\max(b_k,0)+\sum_{j=1}^{k-1}b_j
\right) 
-2\sum_{j=1}^{i-1}b_j,\text{ and}\\*
\mu([\alpha_i]) &= 
\begin{cases}
(-1)^i a_{\lceil i/2 \rceil} + \frac{\mu([\beta_{\lceil
        i/2\rceil}])}{2} & \text{ if $b_{\lceil i/2\rceil}\ge 0$} \\
(-1)^i a_{\lceil i/2 \rceil} + \frac{\mu([\beta_{1+\lceil
        i/2\rceil}])}{2} & 
\text{ if $b_{\lceil i/2\rceil}\le 0$}.
\end{cases}
\end{align*}
Here $\lceil x\rceil$ denotes the smallest integer which is not less than~$x$.
\end{lem}

The proof of this lemma is straightforward. Observe that
if~$\mu([\beta_1])$ is known, then all of the $\mu([\beta_i])$ can be
calculated immediately from the coordinates~$b_j$, and the
$\mu([\alpha_i])$ can then be deduced using the
coordinates~$a_j$. Finally, $\mu([\beta_1])$ can be determined by
using the conditions: that $\mu([\beta_i])\ge 0$ for $1\le i\le n-1$;
that $\mu([\alpha_i]) \ge \left| b_{\lceil i/2\rceil} \right|$ for
$1\le i\le 2n-4$; and that at least one of these inequalities is an
equality (otherwise the foliation would have a leaf parallel to the
boundary of~$D_n$). These conditions give
\[
\mu([\beta_1]) = 2\max_{1\le k\le n-2}\left(
|a_k|+\max(b_k,0)+\sum_{j=1}^{k-1}b_j
\right)
\]
as in the statement of the lemma.

Projectivizing the Dynnikov coordinates yields an explicit
homeomorphism between $S^{2n-5}=(\cC_n)/\R^+$ and the boundary of the
Teichm\"uller space of~$D_n$ (that is, the space of projective
measured foliations on~$D_n$ up to isotopy and Whitehead equivalence).

\begin{remark}
\label{rem:curves}
Let~$\cS_n$ be the set of non-empty unions of pairwise disjoint (but
not necessarily pairwise non-homotopic) essential simple closed curves
on~$D_n$, up to isotopy.  Denote by $S([\alpha])$ the minimum
intersection number of $S\in\cS_n$ with an arc $\alpha\in\cA_n$. Then
there is a bijection $\rho:\cS_n\to\Z^{2n-4}\setminus\{0\}$ defined by
\[\rho(S) = (a,b)=(a_1,\ldots,a_{n-2},\,b_1,\ldots,b_{n-2}),\]
where for $1\le i\le n-2$
\[a_i=\frac{S([\alpha_{2i}])-S([\alpha_{2i-1}])}{2}\qquad
\text{and} \qquad b_i=\frac{S([\beta_i])-S([\beta_{i+1}])}{2}.\]
This bijection is just the restriction of the Dynnikov coordinate
function to the rational measured foliations represented by elements
of~$\cS_n$. 
\end{remark}

\subsection{Update rules}
The Mapping Class Group of~$D_n$ is canonically isomorphic to Artin's
braid group~$B_n$ modulo its center. $B_n$ thus acts on $\fF_n$, and
hence on the space of Dynnikov coordinates. Given~$\beta\in B_n$,
define $\beta:\cC_n\to\cC_n$ by
$\beta(a,b)=\rho\circ\beta\circ\rho\I(a,b)$. 

\begin{remark}
The convention used here for the Artin generators is the normal one in
dynamics, i.e. that used in Birman's book~\cite{birman},
where~$\sigma_i$ denotes the counter-clockwise interchange of the
$i^{\text{th}}$ and $i+1{^\text{th}}$ punctures. Note also the
unfortunate convention that composition is from left to right when
composing braid actions: that is, if $(a,b)\in\cC_n$ and
$\beta_1,\beta_2\in B_n$, then
$(\beta_1\beta_2)(a,b) = \beta_2(\beta_1(a,b))$.
\end{remark}

The {\em update rules} describe the action of the Artin generators
(and their inverses) on~$\cC_n$. For computational and notational
convenience, it is helpful to work in the {\em max-plus semiring}
$(\R,\max,+)$, in which the additive and multiplicative operations are
given by $a\oplus b = \max(a,b)$ and $a\otimes b = a+b$. To simplify
the notation further, formulae in this semiring will use the normal
notation of addition, multiplication, and division, and the fact that
these operations are to be interpreted in their max-plus sense will be
indicated by enclosing the formulae in square brackets. That is,
$[a+b]=\max(a,b)$, $[ab]=a+b$, $[a/b] = a-b$, and $[1]=0$, the
multiplicative identity. For example, the
formula 
\[a_i'=\tropical{\frac{a_{i-1}a_ib_i}{a_{i-1}(1+b_i)+a_i}}\]
given below is just another way of writing
\[a_i'=a_{i-1}+a_i+b_i-\max(a_{i-1}+\max(0,b_i),a_i).\]

\begin{lem}[Update rules for Artin generators]
\label{lem:update}
Let~$(a,b)\in\cC_n$ and \mbox{$1\le i\le n-1$}, and write
$\sigma_i(a,b)=(a',b')$. Then $a_j'=a_j$ and $b_j'=b_j$ except when
$j=i-1$ or $j=i$, and:

\vspace{3mm}

\noindent{if $i=1$ then}
\begin{align*}
a_1' &= \tropical{ \frac{a_1b_1}{a_1+1+b_1} }, &
b_1' &= \tropical{ \frac{1+b_1}{a_1} }; \\[3mm]
\intertext{if $2\le i \le n-2$ then}
a_{i-1}' &= \tropical{ a_{i-1}(1+b_{i-1})+a_ib_{i-1}  }, &
b_{i-1}' &= \tropical{ \frac{a_ib_{i-1}b_i}{a_{i-1}(1+b_{i-1})(1+b_i)
+ a_ib_{i-1}}  }, \\*
a_i' &= \tropical{ \frac{a_{i-1}a_ib_i}{a_{i-1}(1+b_i)+a_i}  }, &
b_i' &= \tropical{ \frac{a_{i-1}(1+b_{i-1})(1+b_i) + a_ib_{i-1}}{a_i} }; \\[3mm]
\intertext{if $i=n-1$ then}
a_{n-2}' &= \tropical{ a_{n-2}(1+b_{n-2})+b_{n-2} }, &
b_{n-2}' &= \tropical{ \frac{b_{n-2}}{a_{n-2}(1+b_{n-2})} }.
\end{align*}
\end{lem}

The update rules for the inverse generators $\sigma_i^{-1}$ can be
obtained from these on conjugating by the involution
\[(a_1,\ldots,a_{n-2},b_1,\ldots,b_{n-2})
\mapsto \tropical{(1/a_1,\ldots 1/a_{n-2},b_1,\ldots,b_{n-2})}\] 
as explained in Section~\ref{sec:updatecontig} below. These rules are
given in the next lemma.

\begin{lem}[Update rules for inverse Artin generators]
\label{lem:updateinverse}
Let~$(a,b)\in\cC_n$ and \mbox{$1\le i\le n-1$}, and write
$\sigma_i\I(a,b)=(a'',b'')$. Then $a_j''=a_j$ and $b_j''=b_j$ except when
$j=i-1$ or $j=i$, and:

\vspace{3mm}

\noindent{if $i=1$ then}
\begin{align*}
a_1'' &= \tropical{ \frac{1+a_1(1+b_1)}{b_1}  }, &
b_1'' &= \tropical{ a_1(1+b_1) }; \\[3mm]
\intertext{if $2\le i \le n-2$ then}
a_{i-1}'' &= \tropical{
  \frac{a_{i-1}a_i}{a_{i-1}b_{i-1}+a_i(1+b_{i-1})}  },  &
b_{i-1}'' &= 
\tropical{ \frac{a_{i-1}b_{i-1}b_i}{a_{i-1}b_{i-1}+a_i(1+b_{i-1})(1+b_i)} }, \\*
a_i'' &= \tropical{ \frac{a_{i-1}+a_i(1+b_i)}{b_i}  }, &
b_i'' &= \tropical{ \frac{a_{i-1}b_{i-1}+a_i(1+b_{i-1})(1+b_i)}{a_{i-1}} };  \\[3mm]
\intertext{if $i=n-1$ then}
a_{n-2}'' &= \tropical{\frac{a_{n-2}}{a_{n-2}b_{n-2}+1+b_{n-2}} }, &
b_{n-2}'' &= \tropical{ \frac{a_{n-2}b_{n-2}}{1+b_{n-2}} }.
\end{align*}
\end{lem}

Using the max-plus notation, the action of any braid~$\beta\in B_n$
on~$\cC_n$ can be computed by composing the functions of
Lemmas~\ref{lem:update} and~\ref{lem:updateinverse} in the normal
way. For a general braid, of course, the resulting rational functions
can be extremely complicated. However, useful results can be obtained for
braids which are ascending or descending sequences of contiguous Artin
generators (or their inverses): these results are described in the
next section.

\section{Update rules for sequences of contiguous generators}
\label{sec:updatecontig}
The update rules for the $n$-braids 
\begin{align*}
\gamma_n^{k,l} &= \sigma_k \sigma_{k+1} \ldots \sigma_{l-1}\sigma_l, \\*
\delta_n^{k,l} &= \sigma_l \sigma_{l-1} \ldots \sigma_{k+1}\sigma_k, \\*
\epsilon_n^{k,l} = \left(\delta_n^{k,l}\right)\I 
&= \sigma_k\I \sigma_{k+1}\I \ldots \sigma_{l-1}\I\sigma_l\I,
\quad\text{and} \\*
\zeta_n^{k,l} = \left(\gamma_n^{k,l}\right)\I  
&= \sigma_l\I \sigma_{l-1}\I \ldots \sigma_{k+1}\I\sigma_k\I,
\end{align*}
where $1\le k\le l\le n-1$, have a relatively simple form. Their
description is, however, complicated by the need to consider
separately the ``end'' cases $k=1$ and $l=n-1$.

\begin{lem}[Update rules for $\gamma_n^{k,l}$]
\label{lem:gamma}
Let~$n\ge 3$, and for $1\le k\le l\le n-1$ let $\gamma^{k,l}_n$
denote the braid $\sigma_k\sigma_{k+1}\ldots \sigma_{l-1}\sigma_l\in B_n$.

Given $(a,b)\in\cC_n$ and an integer~$j$ with $k-1\le j\le n-2$,
write
\[P_j = P_j(b,k) = \tropical{(1+b_{k-1})\prod_{i=k}^j b_i}.\]
(Note the interpretation of this formula in special cases:
$P_j(b,k)=\tropical{\prod_{i=k}^j b_i}$ if~$k=1$, 
$P_j(b,k)=\tropical{(1+b_{k-1})}$ if $j=k-1$, and
$P_j(b,k)=\tropical{1}$ if $k=1$ and $j=0$.)  Similarly, for $k\le
j\le n-2$, write
\[
S_j = S_j(a,b,k) = \tropical{\sum_{i=k}^j \frac{(1+b_i)P_{i-1}}{a_i}}.
\]

Let $(a',b') = \gamma^{k,l}_n(a,b)$. Then $a_j'=a_j$ and $b_j' = b_j$
for $j<k-1$ and for $j > l$.  Moreover,

\vspace{3mm}

\noindent\textbf{1.} If $k>1$ and $l<n-1$ then
\begin{align*}
a_{k-1}' &= \tropical{
a_{k-1}(1+b_{k-1})+a_kb_{k-1}
}, &
b_{k-1}' &= \tropical{
\frac{a_kb_{k-1}b_k}{a_{k-1}(1+b_{k-1})(1+b_k)+a_kb_{k-1}}
}, \\*
a_j' &= \tropical{
a_{j+1}b_{k-1} + a_{k-1}(a_{j+1}S_j + P_j)
},  &
b_j' &= \tropical{
b_{j+1}\left(
\frac
{b_{k-1}+a_{k-1}S_j}
{b_{k-1}+a_{k-1}S_{j+1}}
\right)
} \quad (k\le j<l), \\*
a_l' &= \tropical{
\frac
{a_{k-1}P_l}
{1+b_{k-1}+a_{k-1}S_l}
}, &
b_l' &= \tropical{
b_{k-1} + a_{k-1}S_l
}. \\[3mm]
\intertext{\textbf{2.} If $k>1$ and $l = n-1$ then the formulae in case~1
  hold for $k-1 \le j < n-2$, while}
a_{n-2}' &= \tropical{
b_{k-1}+a_{k-1}\left(S_{n-2}+P_{n-2}\right)
}, &
b_{n-2}' &= \tropical{
\frac{1}{P_{n-2}}\left(
\frac{b_{k-1}}{a_{k-1}} + S_{n-2}
\right)
}. \\[3mm]
\intertext{\textbf{3.} If $k=1$ and $l<n-1$ then}
a_j' &= \tropical{
P_j + a_{j+1}S_j
}, &
b_j' &= \tropical{
b_{j+1}S_j/S_{j+1}
} \qquad (1\le j < l), \\*
a_l' &= \tropical{
P_l/(1+S_l)
},  &
b_l' &= \tropical{
S_l
}. \\[3mm]
\intertext{\textbf{4.} If $k=1$ and $l=n-1$ then the formulae in case~3
  hold for $1\le j < n-2$, while}
a_{n-2}' &= \tropical{
P_{n-2}+S_{n-2}
}, &
b_{n-2}' &= \tropical{
S_{n-2}/P_{n-2}
} .
\end{align*}
\end{lem}
\begin{proof}
The proof is a straightforward induction on~$l\ge k$ for each~$k$,
with the base case~$l=k$ given by the update rules for single braid
generators (Lemma~\ref{lem:update}).

Take, for example, $1<k<n-1$ (cases~1 and~2). Putting~$l=k$ gives
\mbox{$P_l=\tropical{(1+b_{k-1})b_k}$} and
$S_l=\tropical{(1+b_{k-1})(1+b_k)/a_k}$. The rules for $a_{k-1}'$ and
$b_{k-1}'$ given in case~1 of the lemma are identical to those of
Lemma~\ref{lem:update}, while
\begin{align*}
a_k' = a_l'
&= \tropical{
\frac
{a_{k-1}P_l}
{1+b_{k-1}+a_{k-1}S_l}
}
=\tropical{\frac{a_{k-1}(1+b_{k-1})b_k}{(1+b_{k-1}+a_{k-1}(1+b_{k-1})(1+b_k)/a_k}}
\\*
&=\tropical{\frac{a_{k-1}a_kb_k}{a_k+a_{k-1}(1+b_k)}},
\quad\text{and}\\
b_k' = b_l'
&=\tropical{b_{k-1}+a_{k-1}(1+b_{k-1})(1+b_k)/a_k} =
\tropical{\frac{a_kb_{k-1}+a_{k-1}(1+b_{k-1})(1+b_k)}{a_k}},
\end{align*}
in agreement with Lemma~\ref{lem:update}.

Now assume the result is true for some~$l$ with $k\le l < n-1$, so
that \mbox{$\delta_n^{k,l}(a,b)=(a',b')$} as given by case~1 of the
lemma. Let~$(a'',b'')=\delta_n^{k,l+1}(a,b)$, so that \mbox{$(a'',b'')
  = \sigma_{l+1}(a',b')$}. In particular, $a''_j=a'_j$ and
$b''_j=b'_j$ for all $j$ except~$l$ and~$l+1$. Consider $a_{l+1}''$
for $l+1<n-1$ and $a_l''$ for $l+1=n-1$: the other coordinates work
similarly. 

If $l+1 < n-1$, then Lemma~\ref{lem:update} gives
\begin{align*}
a_{l+1}'' &= \tropical {
\frac{a_l'a_{l+1}'b_{l+1}'}
{a_l'(1+b_{l+1}') + a_{l+1}'}
} \\*
&= \tropical {
\frac{a_{l+1}b_{l+1}a_{k-1}P_l/(1+b_{k-1}+a_{k-1}S_l)}
{a_{l+1}+(1+b_{l+1})a_{k-1}P_l/(1+b_{k-1}+a_{k-1}S_l)}
} \\*
&= \tropical {
\frac{a_{k-1}P_{l+1}}
{1+b_{k-1}+a_{k-1}(S_l + (1+b_{l+1})P_l/a_{l+1})}
} \\*
&= \tropical {
\frac{a_{k-1}P_{l+1}}
{1+b_{k-1}+a_{k-1}S_{l+1}}
}
\end{align*}
as required. Similarly if $l+1=n-1$, then Lemma~\ref{lem:update} gives
\begin{align*}
a_l'' &= \tropical {
a_l'(1+b_l')+b_l'} = 
\tropical {
\frac{a_{k-1}P_l}{1+b_{k-1}+a_{k-1}S_l}(1+b_{k-1}+a_{k-1}S_l) + b_{k-1}+a_{k-1}S_l
} \\*
&= \tropical {
b_{k-1}+a_{k-1}(S_{n-2}+P_{n-2})
} \quad\text{as required.}
\end{align*}
\end{proof}

The update rules for $\delta_n^{k,l}$, $\epsilon_n^{k,l}$, and
$\zeta_n^{k,l}$, can be derived from Lemma~\ref{lem:gamma} by
symmetry, conjugating by an appropriate transformation as described
below:
\begin{description}
\item[Reflection in the horizontal diameter of the disk] sends each
  braid generator~$\sigma_i$ to $\sigma_i^{-1}$. The corresponding
  transformation of Dynnikov coordinates is given by
\[(a_1,\ldots,a_{n-2},b_1,\ldots,b_{n-2})\mapsto(-a_1,\ldots,-a_{n-2},b_1,\ldots,
b_{n-2}),\]
or, in max-plus notation,
\[(a_1,\ldots,a_{n-2},b_1,\ldots,b_{n-2}) \mapsto
\tropical{(1/a_1,\ldots,1/a_{n-2},b_1,\ldots,b_{n-2})}.\] Thus the
update rules for~$\epsilon_n^{k,l}$ can be obtained by conjugating the
rules of Lemma~\ref{lem:gamma} by this involution.

\item[Reflection in the vertical diameter of the disk] sends each
  braid generator~$\sigma_i$ to $\sigma_{n-i}^{-1}$. The corresponding
  transformation of Dynnikov coordinates is given by
\[(a_1,\ldots,a_{n-2},b_1,\ldots,b_{n-2})\mapsto(a_{n-2},\ldots,a_1,-b_{n-2},\ldots,
-b_1),\]
or, in max-plus notation,
\[(a_1,\ldots,a_{n-2},b_1,\ldots,b_{n-2}) \mapsto
\tropical{(a_{n-2},\ldots,a_1,1/b_{n-2},\ldots,1/b_1)}.\] Thus the
update rules for~$\zeta_n^{k,l}$ can be obtained by conjugating the
rules of Lemma~\ref{lem:gamma} by this involution.

\item[Rotation through~$\pi$ about the center of the disk] sends each
  braid generator~$\sigma_i$ to $\sigma_{n-i}$. The corresponding
  transformation of Dynnikov coordinates is given by
\[(a_1,\ldots,a_{n-2},b_1,\ldots,b_{n-2})\mapsto(-a_{n-2},\ldots,-a_1,-b_{n-2},\ldots,
-b_1),\]
or, in max-plus notation,
\[(a_1,\ldots,a_{n-2},b_1,\ldots,b_{n-2}) \mapsto
\tropical{(1/a_{n-2},\ldots,1/a_1,1/b_{n-2},\ldots,1/b_1)}.\] Thus the
update rules for~$\delta_n^{k,l}$ can be obtained by conjugating the rules of
Lemma~\ref{lem:gamma} by this involution.
\end{description}

An example which will be used later is given: here the update rules
for $\delta_n^{k,l}$ are derived from those of Lemma~\ref{lem:gamma}
for $\gamma_n^{n-l,n-k}$ by conjugating by a rotation through~$\pi$
about the center of the disk.

\begin{lem}[Update rules for $\delta_n^{k,l}$]
\label{lem:delta}
Let~$n\ge 3$, and for $1\le k\le l\le n-1$ let $\delta^{k,l}_n$ denote
the braid $\sigma_l\sigma_{l-1}\ldots \sigma_{k+1}\sigma_k\in B_n$.

Given $(a,b)\in\cC_n$ and an integer~$j$ with $\max(k-1,1)\le j\le l$
write
\[\wP_j = \wP_j(b,l) = \tropical{(1+b_l)\prod_{i=j}^l\frac{1}{b_i}}.\]
(In the special case $l=n-1$,
$\wP_j(b,n-1)=\tropical{\prod_{i=j}^{n-2}\frac{1}{b_i}}$ for $j<l$,
while \mbox{$\wP_{n-1}(b,n-1)=\tropical{1}$}.)  Similarly, for
$\max(k-1,1)\le j\le l-1$ write
\[
\wS_j = \wS_j(a,b,l) = \tropical{\sum_{i=j}^{l-1}
  \frac{a_i(1+b_i)\wP_{i+1}}{b_i}}.
\]

Let $(a',b') = \delta^{k,l}_n(a,b)$. Then $a_j'=a_j$ and $b_j' = b_j$
for $j<k-1$ and for $j > l$.  Moreover,

\vspace{3mm}

\noindent\textbf{1.} If $k>1$ and $l<n-1$ then
\begin{align*}
a_{k-1}' &= \tropical{
\frac{a_l(1+b_l)+b_l\wS_{k-1}}{b_l\wP_{k-1}}
}, &
b_{k-1}' &= \tropical{
\frac{a_lb_l}{a_l+b_l\wS_{k-1}}
}, \\*
a_j' &= \tropical{
\frac{a_{j-1}a_lb_l}
{a_l+b_l(\wS_j+a_{j-1}\wP_j)}
},  &
b_j' &= \tropical{
b_{j-1}\left(
\frac
{a_l+b_l\wS_{j-1}}
{a_l+b_l\wS_j}
\right)
} \quad (k\le j<l), \\*
a_l' &= \tropical{
\frac
{a_{l-1}a_lb_l}
{a_{l-1}(1+b_l)+a_l}
}, &
b_l' &= \tropical{
\frac
{a_{l-1}(1+b_{l-1})(1+b_l)+a_lb_{l-1}}
{a_l}
}. \\[3mm]
\intertext{\textbf{2.} If $k=1$ and $l < n-1$ then the formulae in case~1
  hold for $2 \le j \le l$, while}
a_1' &= \tropical{
\frac
{a_lb_l}
{a_l+b_l(\wS_1+\wP_1)}
}, &
b_1' &= \tropical{
\frac
{b_l\wP_1}
{a_l+b_l\wS_1}
}. \\[3mm]
\intertext{\textbf{3.} If $k>1$ and $l=n-1$ then}
a_j' &= \tropical{
\frac
{a_{j-1}}
{a_{j-1}\wP_j+\wS_j}
}, &
b_j' &= \tropical{
  b_{j-1}\wS_{j-1}/\wS_j
} \qquad (k\le j \le n-2), \\*
a_{k-1}' &= \tropical{
(1+\wS_{k-1})/\wP_{k-1}
},  &
b_{k-1}' &= \tropical{
1/\wS_{k-1}
}. \\[3mm]
\intertext{\textbf{4.} If $k=1$ and $l=n-1$ then the formulae in case~3
  hold for $2\le j \le n-2$, while}
a_1' &= \tropical{
1/(\wP_1+\wS_1)
}, &
b_1' &= \tropical{
\wP_1/\wS_1
} .
\end{align*}
\end{lem}

\section{Computing topological entropy in families of braids}
\label{sec:main}

If~$\beta\in B_n$ is a pseudo-Anosov braid, then there is some
$(a^u,b^u)\in\cC_n$ (corresponding to the unstable foliation
of~$\beta$) and a number $r>1$ (the dilatation of~$\beta$) such that
$\beta(a^u,b^u)=r(a^u,b^u)$. In this case~$\beta$ has topological
entropy~$h(\beta)=\log r$; there is an element $(a^s,b^s)$ of $\cC_n$
(corresponding to the stable foliation of~$\beta$) with
$\beta(a^s,b^s)=\frac{1}{r}(a^s,b^s)$; and any $(a,b)\in\cC_n$
satisfying $\beta(a,b)=k(a,b)$ for some~$k>0$ is a multiple either of
$(a^u,b^u)$ or of $(a^s,b^s)$.

$\beta$ is a reducible braid if and only if there is some
$(a,b)\in\Z^{2n-4}\setminus\{0\}$ (corresponding to a system of
reducing curves, see Remark~\ref{rem:curves}) with $\beta(a,b)=(a,b)$.

If there is no $(a,b)\in\cC_n$ and~$k>0$ with $\beta(a,b)=k(a,b)$,
then~$\beta$ is a finite order braid, and hence there is some~$N>0$
such that $\beta^N(a,b)=(a,b)$ for all~$(a,b)\in\cC_n$.

In many cases it is possible to do a simultaneous analysis of this
type of every braid in a family. This provides a method of computing
the topological entropy of braids in such families which is more
direct and tractable than the train track approach. In this section,
this method is illustrated with two families of braids considered
in~\cite{kin}, which are of interest in the study of braids of low
topological entropy. These families are $\{\beta_{m,n}\,:\,m,n\ge
1\}$, and $\{\sigma_{m,n}\,:\,1\le m\le n\}$, where
\begin{align*}
\beta_{m,n} &= \sigma_1\ldots \sigma_m \sigma_{m+1}\I\ldots
\sigma_{m+n}\I  = \gamma_{m+n+1}^{1,m}\epsilon_{m+n+1}^{m+1,m+n}\in
B_{m+n+1}, \quad\text{ and}\\*
\sigma_{m,n} &= \sigma_1\ldots \sigma_m \, \sigma_m\ldots\sigma_1 \,
\sigma_1\ldots \sigma_{m+n} =
\gamma_{m+n+1}^{1,m}\delta_{m+n+1}^{1,m}\gamma_{m+n+1}^{1,m+n} \in B_{m+n+1}.
\end{align*}
The approach taken here can be contrasted with the method of proof of
the same results in~\cite{kin}.

\subsection{A family of pseudo-Anosov braids}
The following result establishes that~$\beta_{m,n}$ is a pseudo-Anosov
braid for all $m,n\ge 1$, and provides a formula for the topological
entropy $h(\beta_{m,n})$.
\begin{thm}[The braids $\beta_{m,n}$]
\label{thm:betamn}
Let $m,n\ge 1$. Then $\beta_{m,n}\in B_{m+n+1}$ is a pseudo-Anosov
braid, whose dilatation~$r$ is the unique root in $(1,\infty)$ of the
polynomial
\[f_{m,n}(r)=(r-1)(r^{m+n+1}-1) - 2r(r^m+r^n).\]
The Dynnikov coordinates $(a,b)\in\cC_{m+n+1}$ of the unstable
invariant measured foliation of~$\beta_{m,n}$ are given by
\begin{align*}
a_i &= 
\begin{cases}
-r(r^n+1)(r^i-1) & \text{ if \ } 1\le
  i\le m-1 \\
-(r^{m+1}-1)(r^{n+1}-1) & \text{ if \  } i =
  m\\
-(r^{m+1}-1)(r^{m+n+1-i}-1)r^{i-m}
  & \text{ if \ } m+1 \le i \le m+n-1,
\end{cases}
\\
b_i &= 
\begin{cases}
-(r-1)(r^n+1)r^{i+1} & \text{ if \ } 1\le
  i\le m-1 \\
-(r+1)(r^{m+1}-1) & \text{ if \ } i =
  m\\
-(r-1)(r^{m+1}-1)r^{i-m} \text{\phantom{Some st\,}}
  & \text{ if \ } m+1 \le i \le m+n-1.
\end{cases}
\end{align*}
\end{thm}

\begin{proof}
$f_{m,n}$ has a root $r>1$ since $f_{m,n}(1)=-4$. It will be shown
  that $\beta_{m,n}(a,b)=r(a,b)$, from which the result (and the
  uniqueness of~$r$) follows. 

 Write $N=m+n+1$ and recall that $\beta_{m,n} =
 \gamma_N^{1,m}\epsilon_N^{m+1,N-1}$. Thus to show that
 $\beta_{m,n}(a,b)=r(a,b)$ it suffices to show that
 $\gamma_N^{1,m}(a,b)=r\delta_N^{m+1,N-1}(a,b)$. It will be shown that
 each side of this equation is equal to $(a',b')$, where
\[(a_j',\,\,b_j') = 
\begin{cases}
(ra_j,\,\,rb_j) \quad& 1\le j < m\\
(a_m+b_m,\,\,r(r^n+1)(r+1)) \quad & j=m\\
(a_j,\,\,b_j) & m < j \le m+n-1.
\end{cases}
\]

Observe that 
\begin{equation}
\label{eq:useful}
ra_{m-1}-a_m+a_1 = f_{m,n}(r)+2r(1+r^n)  = 2r(1+r^n) > 0.
\end{equation}

Consider first $(a',b')=\gamma_N^{1,m}(a,b)$, which is given by
Lemma~\ref{lem:gamma}. The first step is to calculate the
quantities~$P_j$ and~$S_j$ from the statement of Lemma~\ref{lem:gamma}
for $1\le j\le m$.

Now $P_j=\sum_{i=1}^j b_i$, giving $P_j=-r^2(r^n+1)(r^j-1)=ra_j$ for
$1\le j < m$; and hence $P_m = P_{m-1}+b_m=ra_{m-1}+b_m$. On the other
hand,
\[S_j = \max_{1\le i\le j}\left( \max(0,b_i) + P_{i-1} - a_i \right) =
\max_{1\le i\le j} \left(ra_{i-1}-a_i\right) \] (setting $a_0=0$),
since $b_i<0$ for all~$i$. Now $ra_{i-1}-a_i = -a_1$ for all $i<m$, so
$S_j=-a_1$ for $1\le j < m$. Finally $S_m = \max(-a_1,ra_{m-1}-a_m) =
ra_{m-1}-a_m$ by~(\ref{eq:useful}).

Let $1\le j\le m-2$. Then (using case~3 of Lemma~\ref{lem:gamma}) 
\begin{align*}
a_j' &= \max(P_j, a_{j+1}+S_j) = \max(ra_j, a_{j+1}-a_1) =
\max(ra_j,ra_j)=ra_j \quad\text{ and }\\*
 b_j' &= b_{j+1}+S_j-S_{j+1} = b_{j+1}=rb_j \qquad\text{ as required.}
\end{align*}

Let $j=m-1$. Then $a_{m-1}' = \max(P_{m-1}, a_m+S_{m-1}) =
\max(ra_{m-1}, a_m - a_1)=ra_{m-1}$ by~(\ref{eq:useful}), and
$b_{m-1}' = b_m+S_{m-1}-S_m = b_m-a_1-(ra_{m-1}-a_m) = rb_{m-1}$ as
required.

Let $j=m$. Then $a_m' = P_m-\max(0,S_m) = ra_{m-1}+b_m-(ra_{m-1}-a_m)
= a_m+b_m$ as required, while $b_m' = S_m = ra_{m-1}-a_m =
2r(1+r^n)-a_1$ by~(\ref{eq:useful}), giving $b_m' = r(r^n+1)(r+1)$ as
required. 

Now let $(a'',b'')=\delta_N^{m+1,N-1}(a,b)$. Showing that $(a'',b'') =
(a',b')/r$, will complete the proof. The argument, using
Lemma~\ref{lem:delta}, is similar to the first part of the
proof. Calculating the quantities $\wP_j$ and $\wS_j$ from the
statement of Lemma~\ref{lem:delta} gives
\begin{align*}
\wP_j &= r^{j-m}(r^{m+1}-1)(r^{m+n-j}-1), &
\wS_j &= -r^n(r-1)(r^{m+1}-1) \qquad (j>m), \\*
\wP_m &= (r^{m+1}-1)(r^n+1), &
\wS_m &= -(r+1)(r^n+1).
\end{align*}

Then, by case~3 of Lemma~\ref{lem:delta},
\begin{align*}
a_m'' &= \max(0,\wS_m)-\wP_m = -\wP_m = (a_m+b_m)/r, \\*
b_m'' &= -\wS_m = (r^n+1)(r+1), \\
a_{m+1}'' &= a_m-\max(a_m+\wP_{m+1},\wS_{m+1}) = a_m-\wS_{m+1} =
a_{m+1}/r, \\*
b_{m+1}'' &= b_m+\wS_m-\wS_{m+1} = b_{m+1}/r+f_{m,n}(r) = b_{m+1}/r,\\
a_j'' &= a_{j-1}-\max(a_{j-1}+\wP_j,\wS_j) = a_{j-1}-\wS_j = a_j/r
\quad (j > m+1), \text{ and}\\*
b_j'' &= b_{j-1}+\wS_{j-1}-\wS_j = b_{j-1} = b_j/r \quad (j>m+1)
\end{align*}
as required.
\end{proof}

\begin{remark}
\label{rem:comps}
 The proof of Theorem~\ref{thm:betamn} is self-contained. However, one
 might ask how the polynomial~$f_{m,n}$ and the Dynnikov coordinates
 of the unstable measured foliation~$[\cF_{m,n},\mu_{m,n}]$
 of~$\beta_{m,n}$ were found.

To find the train tracks for an infinite family of braids, the usual
method would be to compute train tracks (using, for example, the
Bestvina-Handel algorithm~\cite{BH}) for enough examples to spot a
general pattern, and then to prove that the conjectured pattern does
indeed hold for all braids in the family. The method here is
similar. Since~$[\cF_{m,n},\mu_{m,n}]$ is an attracting fixed point
for the action of $\beta_{m,n}$ on the boundary of Teichm\"uller
space, it is easy to find its Dynnikov coordinates numerically. Having
done this for several cases of~$m$ and~$n$, one can guess how the
various maxima in the statements of Lemma~\ref{lem:gamma} and
Lemma~\ref{lem:delta} are resolved. This yields the following
statement (provided $m,n\ge2$):

Assume that $a_i\le 0$; $b_i\le 0$; $a_{i+1}= a_i+b_i$ for $1\le i\le
m-2$; $a_m \le a_{m-1}+b_{m-1}$; $a_{m+1} \ge a_m+b_m$; $a_{i+1} =
a_i-b_i$ for $m+1\le i\le m+n-2$; and $a_{m+n-1} \le b_{m+n-1}$.

Let~$\xi=-a_1+(a_{m-1}+b_{m-1}-a_m)+(a_{m+1}-a_m-b_m)\ge 0$. Then
\[\beta_{m,n}(a,b)=(a',b'),\] where
\begin{align*}
a_i' &=
\begin{cases}
b_1 & i=1\\*
a_{i+1}-a_1 & 2\le i\le m-2\\*
a_{m-1}+b_{m-1}-a_1 & i=m-1\\*
a_{i+1}-\xi & m\le i\le m+n-2\\*
a_{m+n-1}-b_{m+n-1}-\xi & i=m+n-1,
\end{cases}
\\
b_i' &=
\begin{cases}
b_{i+1} & 1\le i\le m-2\\*
a_m-a_{m-1}+b_m-b_{m-1} & i=m-1\\*
a_m-a_{m+1}+b_m+b_{m+1} & i=m\\*
b_{i+1} & m+1\le i\le m+n-2\\*
a_{m+n-1}-b_{m+n-1} & i=m+n-1.
\end{cases}
\end{align*}
Solving for an eigenvalue $r\in(1,\infty)$ and the associated
eigenvector~$(a,b)$ yields the statement of Theorem~\ref{thm:betamn}.
\end{remark}

\begin{remark}
The singularity structure of the invariant foliation
$[\cF_{m,n},\mu_{m,n}]$ can be seen in its Dynnikov coordinates. The
equations
\begin{align*}
a_{i+1} &= a_i+b_i \qquad \text{ if } 1\le i\le m-2 \\*
a_{i+1} &= a_i-b_i \qquad \text{ if } m+1\le i\le m+n-2
\end{align*}
of Remark~\ref{rem:comps} correspond to the existence of an
$(m+1)$-pronged singularity and an $(n+1)$-pronged singularity
respectively.
\end{remark}

\subsection{The reducible case}
Consider now the braids $\sigma_{m,n}$ for $1\le m\le n$. If $n\ge
m+2$ then $\sigma_{m,n}$ is a pseudo-Anosov braid: the following
result can be proved analogously to Theorem~\ref{thm:betamn}.
\begin{thm}[The braids $\sigma_{m,n}$ for $n\ge m+2$]
\label{thm:sigmamn}
Let $1\le m\le n-2$. Then $\sigma_{m,n}\in B_{m+n+1}$ is a pseudo-Anosov
braid, whose dilatation~$r$ is the unique root in $(1,\infty)$ of the
polynomial
\[g_{m,n}(r)=(r-1)(r^{m+n+1}+1) + 2r(r^m-r^n).\]
The Dynnikov coordinates $(a,b)\in\cC_{m+n+1}$ of the unstable
invariant measured foliation of~$\sigma_{m,n}$ are given by
\begin{align*}
a_i &= 
\begin{cases}
r(r^n-1)(r^{i+1}-1) & \text{ if \ } 1\le i\le m-1\\*
(r^{m+1}-1)(r^{m+n-i}-1)r^{i+1-m} & \text{ if \ } m\le i \le m+n-1,
\end{cases}
\\
b_i &= 
\begin{cases}
(r-1)(r^n-1)r^{i+1} & \text{ if \ }1\le i\le m-1\\*
(r-1)(r^{m+1}-1)r^{i-m} \text{\phantom{Some st\,\,}}& \text{ if \ }m\le i\le m+n-1.
\end{cases}
\end{align*}
\end{thm}

However, the focus in this subsection is on the case $n=m+1$, when
$\sigma_{m,n}$ is a reducible braid. Again, the emphasis in the next
result is on the transparent computational nature of the proof, when
compared with a more direct approach such as conjugating the braids in
some suitable way and then appealing to the reader to observe that the
resulting braids leave a certain system of curves invariant.

\begin{figure}[htbp]
\begin{center}
\pichere{0.45}{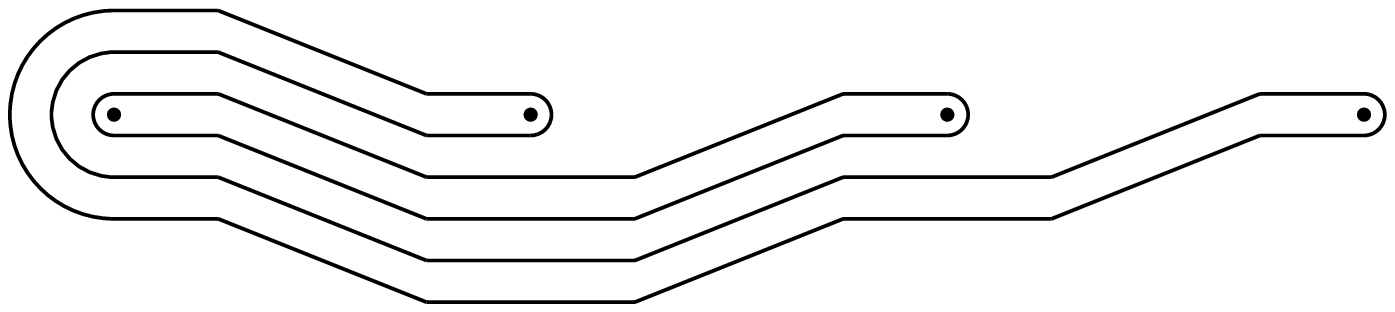}\quad\pichere{1.1}{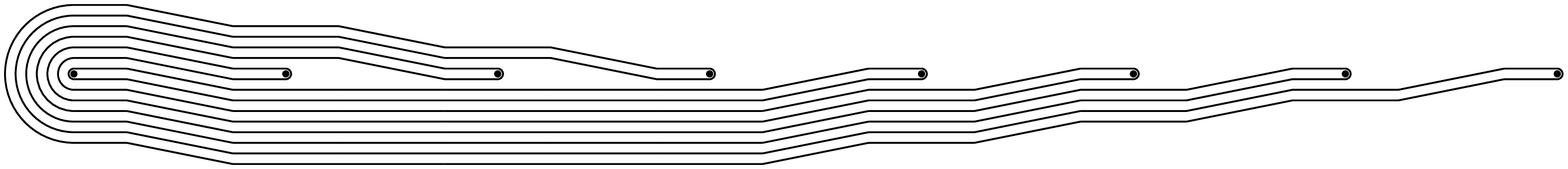}
\end{center}
\caption{The reducing systems $S_1\in\cS_4$ and $S_3\in\cS_8$.}
\label{figure:red}
\end{figure}

\begin{thm}
\label{thm:reducible}
Let $m\ge 1$. Then the braid $\sigma_{m,m+1}\in B_{2m+2}$ is
reducible, having a system of reducing curves $S_m\in\cS_{2m+2}$ with
\mbox{$\rho(S_m)=(a,b)\in \Z^{4m}\setminus\{0\}$} given by
\[
(a_i,b_i) = 
\begin{cases}
(i+1,\,1) & 1\le i\le m \\*
(2m+1-i,1) & m+1\le i\le 2m
\end{cases}
\]
(see Figure~\ref{figure:red}).
\end{thm}

\begin{proof}
Recall that $\sigma_{m,m+1} =
\gamma_{2m+2}^{1,m}\delta_{2m+2}^{1,m}\gamma_{2m+2}^{1,2m+1}$. The
method of proof is to compute successively
$(\A1,\B1)=\gamma_{2m+2}^{1,m}(a,b)$,
$(\A2,\B2)=\delta_{2m+2}^{1,m}(\A1,\B1)$, and
$(\A3,\B3)=\gamma_{2m+2}^{1,2m+1}(\A2,\B2)$, and then to observe that
$(\A3,\B3)=(a,b)$. The calculations are straightforward using
Lemmas~\ref{lem:gamma} and~\ref{lem:delta}.

\noindent\textbf{1.\quad} $(\A1,\B1)$ is computed using case~3 of
Lemma~\ref{lem:gamma}. The quantities $P_j$ and~$S_j$ are given for
$j\le m$ by $P_j=\sum_{i=1}^j b_i=j$ and 
\[S_j = \max_{1\le i\le j}\left(\max(b_i,0)+P_{i-1}-a_i)\right) = 
\max_{1\le i\le j}(1+(i-1)-(i+1)) = -1.\]
Then for $1\le j<m$
\begin{align*}
\A1_j &= \max(P_j,a_{j+1}+S_j) = \max(j, j+2-1) = j+1, \\*
\B1_j &= b_{j+1}+S_j-S_{j+1} = 1-1+1 = 1.
\end{align*}
Finally $\A1_m = P_m-\max(S_m,0)=m-\max(-1,0) = m$, and $\B1_m = S_m =
-1$. Thus
\[(\A1_i,\B1_i) = 
\begin{cases}
(i+1,\,1) & 1\le i < m\\*
(m, -1) & i=m\\*
(2m+1-i,1) & m+1\le i\le 2m.
\end{cases}
\]

\noindent\textbf{2.\quad} $(\A2,\B2)$ is computed using case~2 of
Lemma~\ref{lem:delta}. The quantities $\wP_j$ and $\wS_j$ are given
for $j\le m$ by 
\[\wP_j = \max(\B1_m,0)-\sum_{i=j}^m\B1_i = 1+j-m
\]
and $\wS_j = \max_{j\le i\le m-1}(\A1_i+\max(\B1_i,0)+\wP_{i+1}-\B1_i)
= m+1$. Hence
\begin{align*}
\A2_1 &= \A1_m+\B1_m-\max(\A1_m,\B1_m+\max(\wS_1,\wP_1)) \\*
&= m-1-\max(m,-1+\max(m+1,2-m)) = -1, \\*
\B2_1 &= \B1_m+\wP_1-\max(\A1_m,\B1_m+\wS_1)\\* 
&= -1+(2-m)-\max(m,-1+m+1) = 1-2m,\\
\A2_m &= \A1_{m-1}+\A1_m+\B1_m-\max(\A1_{m-1}+\max(\B1_m,0), \A1_m)\\*
&= m + m -1 -\max(m,m) = m-1,\\*
\B2_m &= \max(\A1_{m-1}+\max(\B1_{m-1},0)+\max(\B1_{m},0),
\A1_m+\B1_{m-1})-\A1_m\\*
&= \max(m+1+0, m+1)-m =  1,\\
\intertext{and for $2\le j<m$}
\A2_j &= \A1_{j-1}+\A1_m+\B1_m-\max(\A1_m,\B1_m+\max(\wS_j,
\A1_{j-1}+\wP_j))\\*
&= j+m-1-\max(m,-1+\max(m+1, 2j+1-m)) = j+m-1-m = j-1,\\*
\B2_j &= \B1_{j-1}+\max(\A1_m,\B1_m+\wS_{j-1}) -
\max(\A1_m,\B1_m+\wS_j) = \B1_{j-1} = 1.
\end{align*}
Thus
\[(\A2_i,\B2_i) = 
\begin{cases}
(-1,\,1-2m) & i=1\\*
(i-1,\,1) & 2\le i \le m\\*
(2m+1-i,1) & m+1\le i\le 2m.
\end{cases}
\]

\noindent\textbf{3.\quad} $(\A3,\B3)$ is computed using case~4 of
Lemma~\ref{lem:gamma}. The quantities $P_j$ and $S_j$ are given by
$P_j=\sum_{i=1}^j\B2_i = j-2m$ (and $P_0=0$); and 
\[S_j=\max_{1\le i\le j}(\max(\B2_i,0)+P_{i-1}-\A2_i).\] 
Now $\max(\B2_i,0)+P_{i-1}-\A2_i$ is equal to~$1$ when $i=1$ and is
negative for $i>1$, and hence $S_j=1$ for all~$j$. Thus $\A3_{2m} =
\max(P_{2m},S_{2m})=1$, $\B3_{2m} = S_{2m}-P_{2m} = 1$, and for $1\le
j<2m$\begin{align*} \A3_j &= \max(P_j,\A2_{j+1}+S_j) = \max(j-2m,
  \A2_{j+1}+1) = \A2_{j+1} + 1 \\* &= \begin{cases} j+1 & 1\le j\le
    m-1 \\* 2m+1-(j+1)+1 & m\le j\le 2m-1,
\end{cases}\\*
\B3_j &= \B2_{j+1}+S_j-S_{j+1} = \B2_{j+1} = 1.
\end{align*}
Hence $(\A3,\B3)=(a,b)$ as required.
\end{proof}

\bibliographystyle{elsart-num-sort}

\end{document}